\documentclass[11pt,reqno]{amsart}
\usepackage{amsmath,amsthm,amsfonts,amssymb,mathrsfs,bm,graphicx,stmaryrd}
\usepackage{mathtools}
\usepackage{dsfont}
\usepackage{multicol}
\usepackage[colorlinks=true,linkcolor=blue]{hyperref}
\usepackage[letterpaper,hmargin=1.0in,vmargin=1.0in]{geometry}
\parskip 	\smallskipamount

\newtheorem{theorem}{Theorem}[section]
\newtheorem{lemma}[theorem]{Lemma}
\newtheorem{corollary}[theorem]{Corollary}
\newtheorem{proposition}[theorem]{Proposition}

\def\N{\mathbb{N}}
\def\P{\mathbb{P}}
\def\Z{\mathbb{Z}}

\def\E{\mathbb{E}}
\newcommand{\cA}{\mathcal{A}}
\newcommand{\cB}{\mathcal{B}}
\newcommand{\cC}{\mathcal{C}}
\newcommand{\cD}{\mathcal{D}}
\newcommand{\cE}{\mathcal{E}}

\newcommand{\cL}{\mathcal{L}}

\newcommand{\bfn}{\mathbf{n}}

\newcommand{\Pc}[1]{\bar{\P}\br{#1}}
\newcommand{\bc}[1]{\left\{{#1}\right\}}
\newcommand{\br}[1]{\left({#1}\right)}
\newcommand{\bs}[1]{\left[{#1}\right]}
\newcommand{\abs}[1]{\left| {#1} \right|}

\begin{document}
\title[Geodesics in LPP]{Sharp deviation bounds for midpoint and endpoint of geodesics in exponential last passage percolation}
\author{Pranay Agarwal}
\address{Pranay Agarwal, Department of Mathematics, University of Toronto, Toronto, Canada.}
\email{pranay.agarwal@mail.utoronto.ca}
\author{Riddhipratim Basu}
\address{Riddhipratim Basu, International Centre for Theoretical Sciences, Tata Institute of Fundamental Research, Bangalore, India} 
\email{rbasu@icts.res.in}

\date{}

\begin{abstract}
    For exponential last passage percolation on the plane we analyze the probability that the point-to-line geodesic exhibits an atypically large transversal fluctuation at the endpoint. We extend the analysis to the probability that the point-to-point geodesic exhibits an atypically large transversal fluctuation at the halfway point. In particular, we show that $p^*_n(t)$, the probability that the point-to-line geodesic from the origin to the line $x+y=2n$ ends at $(n-t(2n)^{2/3}, n+t(2n)^{2/3})$ is such that $n^{2/3}p^*_n(t)=\exp(-(\frac{4}{3}+o(1))t^{3})$ for $t$ large and $p_{n,\frac{1}{2}}(t)$, the probability that the geodesic from the origin to the point $(n,n)$ passes through the point $(\frac{1}{2}n-tn^{2/3}, \frac{1}{2} n+tn^{2/3})$, satisfies $n^{2/3}p_{n,\frac{1}{2}}(t)=\exp(-(\frac{8}{3}+o(1))t^3)$ for $t$ large. The latter result solves a special case of a conjecture from \cite{Liu22}. 
\end{abstract}

\maketitle

\section{Introduction}
We consider the exponential last passage percolation (LPP) model on $\Z^2$, a canonical exactly solvable model in the (1+1)-dimensional KPZ universality class. In exponential LPP one puts independent and identically distributed (i.i.d.)\ exponential weights on the vertices of $\Z^2$, and considers the last passage time between two (ordered) vertices defined by the maximum weight of an up/right path joining the two vertices. Geodesics (i.e., paths attaining the aforementioned maximum) are important objects to study the model and it is known that the fluctuations of the geodesics are governed by the \emph{transversal fluctuation exponent} $2/3$, i.e., the geodesic between say $\mathbf{0}$ and $\mathbf{n}$ ($\mathbf{r}$ shall denote the vertex $(r,r)$ in $\Z^2$) is typically at most $O(n^{2/3})$ distance away from the straight line joining its endpoints. With appropriate scaling, the passage time field converges to a random field called the \emph{directed landscape} \cite{DOV18} and it is known \cite{DV21} that the geodesics in exponential LPP converges to geodesics in the directed landscape. 

In this paper, our aim is to establish sharp bounds for atypically large (at the transversal fluctuation scale) deviations of the point to point geodesic halfway through its journey. That is, for the geodesic $\Gamma_{n}$ from $\mathbf{0}$ to $\mathbf{n}$, let $\Gamma_n(n/2)=(v_1,v_2)$ denote its intersection with the line $\cL_{n/2}:=\{x+y=n\}$. Set $|\Gamma_n(n/2)|:=|v_1-v_2|$. It is known \cite{BSS14, BGZ21, HS18, Aga23, Shen24} that 
$$\exp(-c_2t^3) \le \P(|\Gamma_n(n/2)|\ge tn^{2/3})\le \exp(-c_1t^3),$$
for some constants $c_1,c_2$ whereas it is also known \cite{BB21} that 
$\P(\Gamma_n(n/2)=v)=\Theta(n^{-2/3})$ for any point $v\in \cL_{n/2}$ within $O(n^{2/3})$ distance of the point $(\frac{n}{2},\frac{n}{2})$ (without loss of generality we shall assume $n$ is even, the reader can easily check that this restriction can be removed by adding floor and ceiling signs suitably throughout our arguments). Let $p_{n,\frac{1}{2}}(t)$ denote the probability that $\Gamma_n(n/2)= v_{n,t}=(\frac{n}{2}-\lfloor tn^{2/3} \rfloor, \frac{n}{2}+\lfloor tn^{2/3} \rfloor)$ (we shall omit the floor signs from now on). Based on the above results, it is natural to expect that $n^{2/3}p_{n,\frac{1}{2}}(t)=\exp(-(c+o(1))t^{3})$ for some $c>0$. Our first theorem establishes this with $c=\frac{8}{3}$.

\begin{theorem}
\label{t:p2p}
Let $\epsilon>0$ be fixed. There exists $t_0, N_0$ such that for $t>t_0$, $n>N_0 \vee t^{40000}$ we have 
$$ n^{-2/3}\exp(-(\frac{8}{3}+\epsilon)t^3)\le p_{n}(t) \le  n^{-2/3}\exp(-(\frac{8}{3}-\epsilon)t^3).$$
\end{theorem}

An easy corollary of the above theorem is the following. 

\begin{corollary}
    \label{c:liu}
    Let $\epsilon>0$ be fixed. There exists $t_0, N_0$ such that for $t>t_0$, $n>N_0 \vee t^{40000}$ we have
    $$\exp(-(\frac{8}{3}+\epsilon)t^3)\le \P(\Gamma_{n}~\text{lies above}~v_{n,t}) \le \exp(-(\frac{8}{3}-\epsilon)t^3).$$
\end{corollary}

We also have a similar result for point-to-line geodesics. Let $v^*_{n,t}$ denote the point $(n-t(2n)^{2/3},n+t(2n)^{2/3})$. Let $p^*_{n}(t)$ denote the probability that the point-to-line geodesic from $\mathbf{0}$ to  the line $\cL_{n}:=\{x+y=2n\}$, denoted $\Gamma^*_n$, ends at $v^{*}_{n,t}$. We have the following analogue of Theorem \ref{t:p2p} in this case. 

\begin{theorem}
\label{t:p2l}
For each $\epsilon>0$, there exists $t_0, N_0$ such that for $t>t_0$, $n>N_0 \vee t^{40000}$ we have 
$$ n^{-2/3}\exp(-(\frac{4}{3}+\epsilon)t^3)\le p^*_{n}(t) \le  n^{-2/3}\exp(-(\frac{4}{3}-\epsilon)t^3).$$
\end{theorem}

As before, we have the following corollary. 

\begin{corollary}
    \label{c:p2l}
    For each $\epsilon>0$, there exists $t_0, N_0$ such that for $t>t_0$, $n>N_0 \vee t^{40000}$ we have
    $$\exp(-(\frac{4}{3}+\epsilon)t^3)\le \P(\Gamma^*_{n}~\text{lies above}~v^*_{n,t}) \le \exp(-(\frac{4}{3}-\epsilon)t^3).$$
\end{corollary}
{We note that the upper bound on the range of $t$ (i.e. $t\le n^{1/40000}$) in all the above results is suboptimal. We have made no effort to optimize this since our interest is in the event of large on-scale transversal fluctuations (i.e. $t$ large but fixed).}

\subsection{Outline of the argument and a discussion of related results}
The proofs of Theorems \ref{t:p2p} and \ref{t:p2l} are similar so we shall only discuss the first case. The main step is to prove that the probability that $\Gamma_n(n/2)$ lies on an interval of length $O(n^{2/3})$ (call this interval $I$) on $\cL_{n/2}$ with center at $v_{n,t}$ is $\exp(-(8/3+o(1))t^3)$, Theorem \ref{t:p2p} then follows from an averaging argument. The upper bound and the lower bound are proved via separate arguments.

The idea for the upper bound is as follows. If $\Gamma_n(n/2)\in I$, it is clear that $\sup_{v\in I} T_{\mathbf{0},v}+T_{v,\mathbf{n}}$ cannot be much smaller than $4n$. Since $\E T_{\mathbf{0},v}=\E T_{v,\mathbf{n}} \approx 4n-\Theta(t^2)n^{1/3}$, it is clear that the probability of the event $\{\Gamma_n(n/2)\in I\}$  can be upper bounded by the probability that $\sup_{v\in I} T_{\mathbf{0},v}+T_{v,\mathbf{n}}$ is atypically large. Since the two summands are independent, one can estimate this probability by using the sharp moderate deviation estimates of passage times (see e.g.\ \cite{BBBK24}) and this leads to the claimed upper bound. 

For the lower bound we construct the following event. We ask that there are two paths, from $\mathbf{0}$ to $v_{n,t}$ and $v_{n,t}$ to $\mathbf{n}$ respectively, each of which have length $\approx 2n$ (i.e., $\Theta(t^2)n^{1/3}$ larger than typical) such that the excess length is approximately uniformly distributed along their journeys. It follows from the sharp moderate deviation estimates that the probability of this event is at least $\exp(-(8/3+o(1))t^3)$. The non-trivial part is to show that this event indeed implies that $\Gamma_n$ approximately follows the concatenation of these planted paths and hence we get the desired lower bound for the probability of the event $\Gamma_n(n/2)\in I$. 

One can actually prove something stronger for the upper bound. For $\gamma\in (0,1)$ fixed, one can consider the probability that the geodesic $\Gamma_n$ passes above the point $v_{n,\gamma,t}=(\gamma n-tn^{2/3},\gamma n+tn^{2/3})$. Using essentially the same argument as above (and solving a slightly different variational problem) one can actually show that the is probability is upper bounded by  $\exp(-(c_{\gamma}+o(1))t^3)$ where $c_{\gamma}=3^{-1}(\gamma(1-\gamma))^{-3/2}$. In fact it was conjectured by Liu \cite[Conjecture 1.5]{Liu22} that this is the correct leading order behavior of this probability. Therefore, our Corollary \ref{c:liu} proved Liu's conjecture for the special case $\gamma=1/2$. 

However, we believe that Liu's conjecture may not be true for $\gamma\ne 1/2$. Indeed, one can check that the planted path that leads to  solution of the variational problem leading to the upper bound described above does not (if $\gamma\ne 1/2$) imply that the planted path approximates the geodesic $\Gamma_n$. To get the correct leading order behaviour for general $\gamma$ one needs to solve a more complicated variational problem. The solution for this variational problem was conjectured in the recent independent work of Das, Dauvergne and Vir\'ag \cite[Conjecture 2.11]{DDV24}. They also prove a version of Corollary \ref{c:liu} for the directed landscape which, as mentioned before, is the scaling limit of exponential LPP (see \cite[Corollary 2.9]{DDV24}). This result could also be recovered taking limits in our Corollary \ref{c:liu} (the pre-factor $32/3$ in \cite[Corollary 2.9]{DDV24} differs from our $8/3$ due to a scaling difference). In fact, it is conjectured in \cite{DDV24} that for $\gamma\ne 1/2$, the planted path that the geodesic approximately follows under the atypical event of large deviation at an intermediate time $\gamma\ne 1/2$ is no longer piecewise linear, but will have two linear pieces connected by a non-linear curve. The work \cite{DDV24} goes far beyond this problem, and establishes a full large deviation principle for the geodesics as well as upper tail large deviation principle for the directed landscape. 

We also point out that taking limit as $n\to \infty$ in Corollary \ref{c:p2l} above, together with the convergence to point-to-line last passage time profile to Airy$_2$ process minus a parabola (see e.g.\ \cite{BF08}), one could also get the tail behavior for the maximizer of the parabolically shifted Airy$_2$ process, which was previously proved in \cite{BLS12,Sch12}.

\subsection*{Acknowledgments}
We thank Duncan Dauvergne and B\'alint Vir\'ag for informing us of their work \cite{DDV24} and for insightful discussions. RB is partially supported by a MATRICS grant (MTR/2021/000093) from SERB, Govt.\ of India, DAE project no.\ RTI4001 via ICTS, and the Infosys Foundation via the Infosys-Chandrasekharan Virtual Centre for Random Geometry of TIFR.

\subsection*{Organisation of the paper}
In Section \ref{s:prelim} we collect some useful estimates which we shall use throughout. Section \ref{s:p2l} is devoted to the proof of Theorem \ref{t:p2l}.  In Section \ref{s:p2p} we prove Theorem \ref{t:p2p} and Corollary \ref{c:liu}. 

\section{Preliminary estimates}
\label{s:prelim}

In this section, we collect a number of preliminary estimates that will be used throughout. We start with setting up some notation. We define $1_f = (-1, 1)$ and the functions $\phi$ and $\psi$ as
 \begin{align*}
 \phi((x,y)) = x+y && \mathrm{and} && \psi((x,y)) = y-x.
\end{align*}
Recall the definition of $\cL_k$ where $k\in \frac{1}{2}\Z$. The geodesic from a point $v$ to the line $\cL_k$ will be denoted by $\Gamma^*_{v,k}$. The notation is simplified to $\Gamma^*_{k}$ when the starting point is the origin. The weight of a path $\gamma$ will be denoted by $\ell(\gamma)$. We denote the intersection point of a path $\Gamma$ with the line $\cL_k$ by $\Gamma(k)$ and $\abs{\Gamma(k)} := \psi(\Gamma(k))$. For a point to line geodesic ${\Gamma}^*$, ending on the line $\cL_{k}$, we shall use $\abs{{\Gamma}^*}$ to denote $\abs{{\Gamma}^*(k)}$ and for the point to point geodesic $\Gamma$ we shall use $\abs{{\Gamma}}$ to denote,
$$\abs{\Gamma} = \sup \bc{\psi(w) : w \in \Gamma}.$$ 


Let us first set-up the basic estimates for the passage times. For $0 \le t \ll n^{2/3}$. recall the point $v=v^*_{n,t}$ from above. Let $T_{\mathbf{0},v}$ denote the passage time from the origin to $v$. We record the following basic estimates for upper tail deviations of last passage times (see, e.g.\cite{BBBK24} where these estimates were recently proved with optimal constants).  

\begin{proposition}[{\cite[Theorem 1.6]{BBBK24}}]
    \label{t:tails}
    Let $\epsilon>0$ be fixed. For $0\le t \ll n^{1/3}$ and $1\le s\ll n^{2/3}$, and for $n$ sufficiently large we have 
    \begin{enumerate}
        \item[(i)]  $\exp(-(\frac{4}{3}+\epsilon)s^{3/2})\le \P(T_{\mathbf{0},v^{*}_{n,t}}\ge 4n-t^22^{4/3}n^{1/3}+s2^{4/3}n^{1/3})\le \exp(-(\frac{4}{3}-\epsilon)s^{3/2})$.
        \item[(ii)]   $\exp(-(\frac{4}{3}+\epsilon)s^{3/2})\le \P(\sup_{v\in \cL_{n}} T_{\mathbf{0},v} \ge 4n+sn^{1/3}) \le \exp(-(\frac{4}{3}-\epsilon)s^{3/2}).$
    \end{enumerate}
\end{proposition}

For the lower tail, following weaker estimate from \cite{LR10} will suffice for our purposes. 

\begin{proposition}[{\cite[Theorem 4]{LR10}}]
    \label{p:ltails}
    There exists $C>0$ such that for $s\in (1,n^{2/3})$
    $$\P\left(\sup_{v\in \cL_{n}} T_{\mathbf{0},v}\le 4n-sn^{1/3}\right)\le \exp(-cs^3).$$
\end{proposition}

We shall also need upgraded versions of these point-to-point and point-to-line estimates where the starting point is a narrow interval instead of a point. For a $\delta>0$ to be chosen later, let us consider the interval 
$$I=I_{\delta}:=\{(x,-x): |x|\le \delta n^{2/3}\}\subset \cL_0.$$
We have the following result from \cite{BGZ21} (a similar result was also established in \cite{BSS14}).

\begin{proposition}[{\cite[Lemma C.5]{BGZ21}}]
 \label{p:inf}
 There exists $c>0$ such that for all $\delta$ small 
 \begin{enumerate}
     \item $\P[\inf_{v\in I} T_{v,\cL_{n}}\le 4n-sn^{1/3}]\le \exp(-cs^{3})$.
     \item  $\P[\inf_{v\in I} T_{v,\mathbf{n}+v}\le 4n-sn^{1/3}]\le \exp(-cs^{3})$.
 \end{enumerate}
\end{proposition}


We also have the following result.

\begin{proposition}
    \label{p:sup}
    Let $\epsilon>0$ be fixed. There exists $\delta_0,t_0,s_0$ such that For $\delta \le \frac{1}{t} \wedge \delta_0$ and for all $t_0\le t\ll n^{1/3}$ and all $s\in (s_0, n^{2/3})$ we have  
    $$\P\left[\sup_{v\in I} T_{v,v+v^{*}_{n,t}}\ge 4n-t^{2}2^{4/3}n^{1/3}  +s2^{4/3}n^{1/3}\right]\le \exp(-(\frac{4}{3}-\epsilon)s^{3}).$$
\end{proposition}

\begin{proof}
    The proof is same as the proof of \cite[Proposition A.5]{BBF22} (which proves the interval to line version of this lemma); we omit the details. 
\end{proof}


\section{Point-to-line geodesics}
\label{s:p2l}

Our proofs of Theorems \ref{t:p2p} and \ref{t:p2l} are similar, with the former being slightly more complicated. In this section, we shall present the proof for Theorem \ref{t:p2l}, and adapt the arguments to establish Theorem \ref{t:p2p} in the next section. 

Notice that with $p_n^*(t)$ as in the statement of Theorem \ref{t:p2l} and $I$ as defined in the previous section, we have
$$2\delta n^{2/3}p_n^*(t)=\E [\#\{v\in I: \Gamma^*_{v,n}~\text{ends at}~v+v^{*}_{n,t}\}].$$

Let us denote the set inside the expectation by $I_*$ and its cardinality by $\abs{I_*}$. Theorem \ref{t:p2l} will follow from the following two estimates. 


\begin{proposition}
    \label{p:p2lub}
    In the set-up of Theorem \ref{t:p2l}, there exists $c_0>0$ sufficiently small such that for $c<c_0$ and with $\delta=\frac{c}{t}$ we have for all $t$ sufficiently large
    $$\E[|I_*|]\le  \exp(-(\frac{4}{3}-\epsilon)t^3).$$
\end{proposition}

\begin{proposition}
    \label{p:p2llb}
    In the set-up of Theorem \ref{t:p2l} we have for some $\delta<1$
    $$\P[|I_*|\ge 1]\ge  \exp(-(\frac{4}{3}+\epsilon)t^3).$$
\end{proposition}

\subsection{Upper bound}
The proof of Proposition \ref{p:p2lub} is divided into two lemmas. 

\begin{lemma}
    \label{l:p2lub}
    In the set-up of Proposition \ref{p:p2lub} we have
    $$\P[|I_*|\ge 1]\le  \exp(-(\frac{4}{3}-\epsilon)t^3).$$
\end{lemma}

\begin{proof}
    We know from Proposition \ref{p:inf} that for $C$ sufficiently large
    $$\P[\inf_{v\in I} T_{v,\cL_{n}}\le 4n-Ctn^{1/3}]\le \exp(-\frac{4}{3}t^{3}).$$ 
    On the other hand we also know by Proposition \ref{p:sup} that for $\delta=\frac{c}{t}$ sufficiently small we have
    $$\P[\sup_{v\in I} T_{v,v+v^*_{n,t}} \ge 4n- t^{2}2^{4/3}n^{1/3} + (t^{2}2^{4/3}n^{1/3}- Ctn^{1/3})] \le \exp(-(\frac{4}{3}-\epsilon)t^3).$$
    It is clear that on the event 
    $$\{\inf_{v\in I}T_{v,\cL_{n}}> 4n-Ctn^{1/3}\} \cap \{\sup_{v\in I} T_{v,v+v^*_{n,t}} <4n-Ctn^{1/3}\},$$
    We have $I_{*}=\emptyset$. This completes the proof of the proposition. 
\end{proof}

Suppose now that $|I_{*}|\ge k$. By uniqueness and ordering of the geodesics, $|I_{*}|\ge k$ implies that there exist $k$ disjoint geodesics between points in $I$ and points in $v^*_{n,t}+I$. The next lemma is then immediate from \cite[Proposition 3.1]{BHS18} (see also  \cite[Proposition 1.6]{BBB23}). 

\begin{lemma}
    \label{l:multirare}
    There exists $k_0,c>0$ such that for $k_0\le k \le n^{0.01}$
    $$\P[|I_*|\ge k]\le \exp(-ck^{1/128}).$$
\end{lemma}

We can now complete the proof of Proposition \ref{p:p2lub}. 

\begin{proof}[Proof of Proposition \ref{p:p2lub}]
Choose $k_*=Ct^{384}$ where $C$ is sufficiently large such that we have, using Lemma \ref{l:multirare} that
$\P[|I_*|\ge k_*]\le e^{-\frac{4}{3}t^3}$. 
Writing 
$$\E|I_{*}|=\sum_{k\ge 1} \P[|I_*|\ge k]=\sum_{k\le k_*} \P[|I_*|\ge k]+ \sum_{k_*\le k \le n^{0.01}} \P[|I_*|\ge k] + \sum_{n^{0.01}\le k \le \delta n^{2/3}} \P[|I_*|\ge k].$$

Observe that, by our choice of $k_*$, together with Lemmas \ref{l:p2lub} and \ref{l:multirare} we get that the first sum above is bounded by $Ct^{384}\exp(-(\frac{4}{3}-\epsilon)t^3)$, the second sum above is upper bounded by $C\exp(-\frac{4}{3}t^3)$ and the third term, using $n \ge t^{40000}$ is also upper bounded by $C\exp(-\frac{4}{3}t^3)$. This completes the proof of the proposition. 
\end{proof}

\subsection{Lower bound}
For the lower bound, we do the following. Let $v_1$ and $v_2$ denote the left and right endpoint of $I$ respectively. Let $\cA_*$ denote the event that the endpoint of $\Gamma^*_{n,v_1}$ is below $v_1+v^*_{n,t}$ and the endpoint of $\Gamma^*_{n,v_2}$ is above $v_2+v^*_{n,t}$. By ordering of the geodesics it is clear that, on $\cA_*$ we have $|I_*|\ge 1$. Therefore Proposition \ref{p:p2llb} follows from the next result. 

\begin{proposition}
    \label{p:p2llbbasic}
    In the set-up of Proposition \ref{p:p2llb}, we have $\P(\cA_*) \ge \exp(-(\frac{4}{3}+\epsilon)t^{3})$. 
\end{proposition}

\begin{figure}
    \centering
    \includegraphics[scale = 0.3]{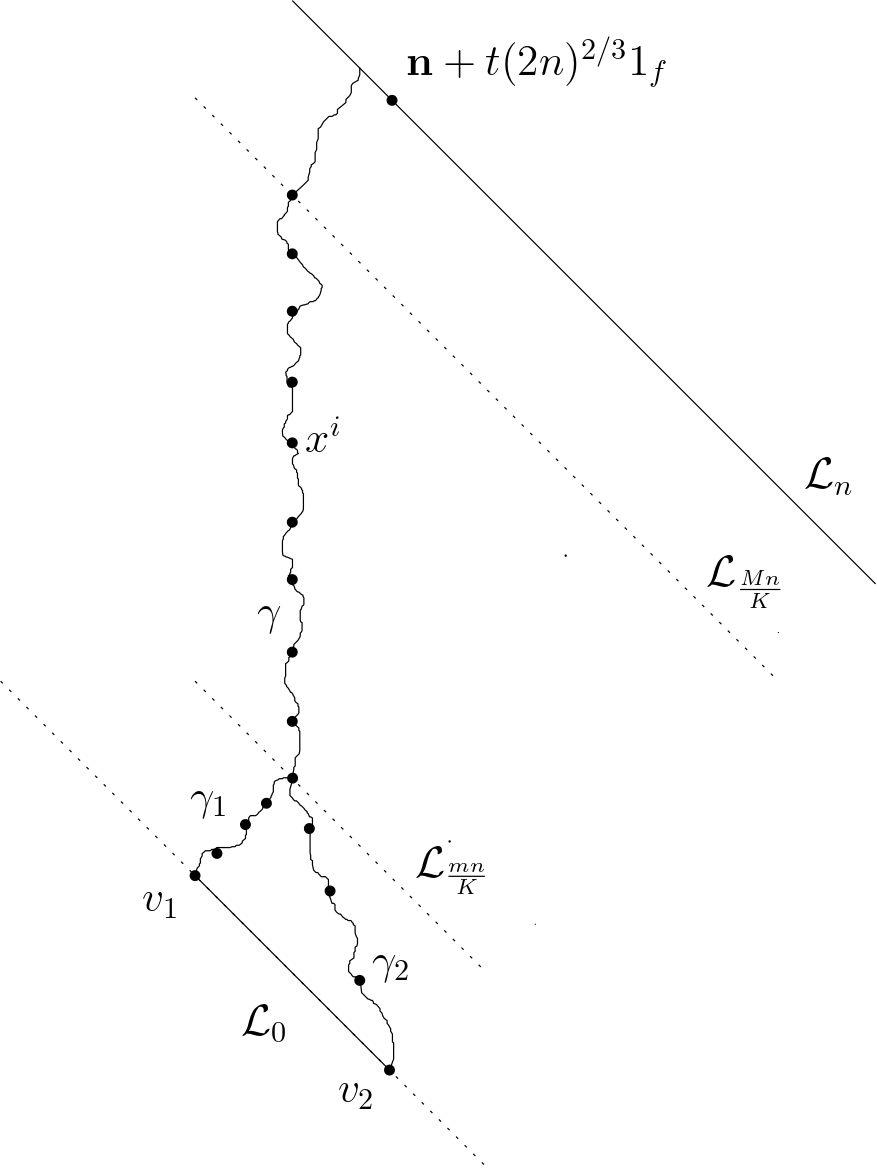}
    \caption{Illustration for Proposition \ref{p:p2llbbasic}. The basic idea here is that we plant two paths, one each from $v_1$ and $v_2$ to $\cL_{n}$ such that the lengths of these paths are larger than $4n$, and that these paths roughly follow the straight-lines joining $v_1$ and $v_2$ to $v^*_{n,t}$. We further ask that the excess weight is approximately uniformly distributed along different parts of these paths. We then show that the probability cost of planting these paths is $\approx\exp(-4t^3/3)$ and that the event $\cA_*$ happens with probability bounded away from $0$ conditional on the planted paths having appropriate weight.}
    \label{fig:pl_const}
\end{figure}

We shall now provide a construction along with a series of lemmas which will together prove Proposition \ref{p:p2llbbasic}. Our construction will depend on a number of parameters chosen as follows. 
\begin{align} \label{para}
    \delta &=  t^{-0.4}, &K=  t^{2.5},\\
    M &= K - t, & m = \frac {M}{2^{2/3} t^{1.4}}.        
\end{align} 
We construct a sequence of points separated by $\frac{n}{K}$ distance in time, going from $v_1$ and $v_2$ to the point $v^*_{n,t}$. {For notational convenience, we shall assume that $K$ and $n/K$ are integers, however it will be clear from the proof that our arguments trivially be adapted to the case when $n/K$ is not an integer.} From the point $v_1$, we have a sequence of points on a line with slope 1, till we intersect the line $\cL_{mn/K}$, at say $x^m$. From $v_2$ we have a sequence of points, which lie on the straight line joining  $v_2$ and $x^m$. Lastly, we have a sequence of points from $x^m$ to the point $\frac{M}{K}\bfn + t(2n)^{2/3}1_f$ along the straight line joining them. We ask that the geodesics connecting consecutive points have more than the typical weight. Let
\begin{align*}
    \cA^1_i &= \bc{T_{x^{i-1}_1,x^{i}_1} \geq 4n/K + \delta t^2 2^{4/3} n^{1/3}/K} &\mathrm{where}& &x^{i}_1 =  v_1 + \frac{i\bfn}{K}, \\
    \cA^2_i &= \bc{T_{x^{i-1}_2,x^{i}_2} \geq 4n/K + \delta t^2 2^{4/3} n^{1/3}/K} &\mathrm{where}&&x^{i}_2 = v_2 + \frac{i \bfn}{K} + \frac{2i}{m}\delta n^{2/3} 1_f,\\
    \cA_i &= \bc{T_{x^{i-1},x^{i}} \geq 4n/K + \delta t^2 2^{4/3} n^{1/3}/K} &\mathrm{where}&&x^{i} = \br{\bfn + \frac{K}{M}t(2n)^{2/3}1_f}\frac{i}{K}.
\end{align*}
We define the events $\cA^1_i$ and $\cA^2_i$ for $0< i \leq m$ and $\cA_i$ for $m < i \leq M$. We complete our construction by defining the event
\[
    \cA_{M+1} = \bc{{T}_{x^M, n}^* \geq (4n + \delta t^2 2^{4/3}n^{1/3})\frac{K-M}{K}},
\]
where ${T}_{x^M, n}^*$ denotes the point to line passage time from $x^{M}$ to $\cL_{n}$. These geodesics concatenate to form paths which can be divided into three parts:
\begin{itemize}
    \item \textbf{Initial Part:} This part consists of two paths $\{\gamma_i\}_{i=1}^{2}$ which start from $v_i$. The two parts eventually meet at the point $x^m$ on the line $\cL_{mn/K}$.
    \item \textbf{Middle Part:} This part is a single path from the intersection point $x_m$ of the initial part to $x^{M}$ on $\cL_{Mn/K}$. 
    \item \textbf{Final Part:} The final part consists of a single point-to-line geodesic from $x^{M}$ to $\cL_{n}$.
\end{itemize}
We refer to the union of all these paths as $\gamma$, which is a path with two branches $\gamma_1$ and $\gamma_2$. See Figure \ref{fig:pl_const} for a visual description of the construction. 

For ease of notation we shall denote the ratio $K/M$ by $\kappa$. Note that $\kappa$ will be very close to $1$ for large $t$. Each parameter in our construction grants us certain control over the construction. The parameter $K$ decides the distance between consecutive points, while $m$ decides the point at which the point sequence from $v_1$ and $v_2$ meet. The parameter $M$ controls the point at which we enter the final part. We want $\delta$ to be small so that the cost of controlling the concatenated paths from $v_1$ and $v_2$ is not large. A large value of $K$ ensures that if $\Gamma^*_n$ is not following the concatenated path, then it has to pay a significant penalty in passage time. We want both the initial and the final parts to be much smaller than the middle part, this will help us ensure that the probability cost of our construction is of the same order as the bound we have set out to prove. 

Let us denote the intersection of all the events which are a part of our construction by $\cA$. We now show that the probability cost for constructing $\gamma$ is of the same order as the desired lower bound. By Proposition \ref{t:tails}, and our choice of the parameters, for arbitrarily small $\epsilon$, we have that for $t$ large enough
\begin{align*}
    \P\br{\cA^1_i} &\ge \exp\br{-\br{\frac{4}{3} +\frac{\epsilon}{2}} \frac{\delta^{3/2}t^3}{K}},
    & \P\br{\cA^2_i} \ge \exp\br{-\br{\frac{4}{3} +\frac{\epsilon}{2}} \frac{t^3 (\delta+4 \kappa^2)^{3/2}}{K}},\\
    \P\br{\cA_i} &\ge \exp\br{-\br{\frac{4}{3} +\frac{\epsilon}{2}} \frac{t^3 (\delta+\kappa ^2)^{3/2}}{K}},
    &\P\br{\cA_{M+1}} \ge \exp \br{-\br{\frac{4}{3} +\frac{\epsilon}{2}} \frac{(K-M)\delta^{3/2} t^3}{K}}.
\end{align*}
Since all of these events are increasing, by the FKG inequality we have that
\begin{align}
\label{e:alb}
    \P\br{\cA} &\geq \exp\br{-c \delta t^2} \exp\br{-\br{\frac{4}{3} +\frac{\epsilon}{2}} (\delta+\kappa ^2)^{3/2} t^3}\nonumber\\
    &\geq \exp\br{-\br{\frac{4}{3} +{\epsilon}} t^3},
\end{align}
by our choice of parameters and taking $t$ large enough.

We shall show that the geodesic $\Gamma^*_{v_i, n}$ for $i=1,2$ coalesce with $\Gamma^*_{x^M, n}$ with large probability given the event $\cA$. Let  $\cB_{v_i}$ denote the event that $\Gamma^*_{v_i,n}$ does not meet with $\Gamma^*_{x^{M},n}$. The following lemma will give us a bound on the conditional probability of $\cB_{v_i}$.

\begin{lemma} \label{l:coal}
    For $i=1,2$ we have
    \[
        \P\br{\cB_{v_i} \:\vert\: \cA} \leq 1/4.
    \]
\end{lemma}

\begin{figure}
    \centering
    \includegraphics[scale = 0.8]{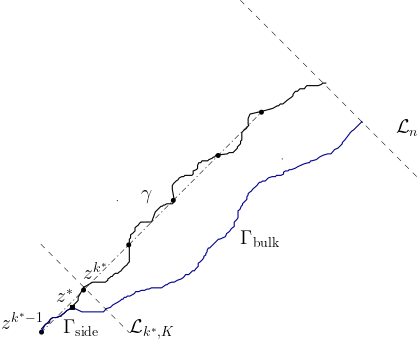}
    \caption{Illustration for Lemma \ref{l:coal}. Let us consider the case that $z^*$ belongs to $\gamma_2$. On the event $\cB_k'$, there exists a path from $z^{k-1}$ to $\cL_{n}$ that does not intersect $\gamma$ after $z^k$ and still has atypically large length. We show using the BK inequality that this event has low probability.}
    \label{fig:p2l_bk}
\end{figure}

\begin{proof}
We shall only show the bound for $\P(\cB_{v_2} \:\vert\: \cA)$ the other one follows similarly. Let $z^*$ be the last intersection point of ${\Gamma}_{v_2,n}^*$ and $\gamma$ and let $k^*$ denote the smallest integer such that $z^*$ lies before $\cL_{\frac{nk^*}{K}} := \cL_{k^*,K}$. On the event $\cB_{v_2}$, $k^*$ must be less than $M$. Thus, for $\cB_k := \{ k^* = k\}$,
    \[
        \cB_{v_2} = \cup_{k \leq M} \cB_k.
    \]  

We further define $\cB'_{k}=\cB_{k}$ if $k\ge m$. If $k<m$ we denote $\cB'_{k}\subset \cB_{k}$ to be the event that the last intersection ${\Gamma}_{v_2,n}^*$ and $\gamma$ belongs to $\gamma_2$. We shall only show  
$\P(\cup_{k \leq M} \cB'_k \:\vert\: \cA)\le \frac{1}{8}$; the remaining case, i.e., showing that $\P(\cup_{k \leq M} \cB_k\setminus \cB'_k \:\vert\: \cA)\le \frac{1}{8}$ is analogous and is sketched without details at the end of the argument. 

 We do a bit of relabeling; set
 \begin{equation*}
    z^i = \begin{cases}
        x^i_2\;\;\;\;\; \mathrm{for} \;\;\;\;\;0 \leq i \leq m,\\
         x^i \;\;\;\;\;  \mathrm{for} \;\;\;\;\;m < i \leq M.
    \end{cases}
\end{equation*}

    On the event, $\cB'_k \cap \cA$, consider the path $\Gamma$ from $z^{k-1}$ to $\cL_n$ which follows $\gamma_2$ until $z^*$ and then follows $\Gamma^*_{v_2,n}$. By definition this path has length larger than $\gamma_2$ from $z^{k-1}$ to $\cL_n$ and hence has a passage time of more than $(4n + {\delta} t^2 2^{1/3} n^{1/3})\frac{K-k+1}{K}$. We divide $\Gamma$ into two parts,
        \begin{align*}
            \Gamma_\mathrm{side} &= \bc{\phi(z^{k-1}) \leq \phi(w) \leq \phi(z^k): w \in {\Gamma}}, \\
            \Gamma_\mathrm{bulk} &= \bc{\phi(z^k) < \phi(w) \leq 2n : w \in {\Gamma}}.
        \end{align*}
    Either $\Gamma_\mathrm{bulk}$ or $\Gamma_\mathrm{side}$ must contribute to the excess weight. Hence, we define the following
        \begin{align*}
            \cB_{k}^1 &= \cB'_{k} \cap \bc{\ell(\Gamma_\mathrm{bulk}) \geq (4n + {\delta} t^2 2^{1/3} n^{1/3})\frac{K-k}{K}}, \\
            \cB_{k}^2 &= \cB'_{k} \cap \bc{\ell(\Gamma_\mathrm{side}) \geq \frac{4n}{K} +{\delta} t^2 2^{4/3} n^{1/3}\frac{K-k}{2K}},\\
            S &= \bc{x \in \cL_{k,K} \:\vert\:\abs{\psi(x-z^{k-1})} \leq n^{2/3}}.
        \end{align*}        
    On the event $\cB_k^1$, either $\Gamma_\mathrm{bulk}$ is a path from $S$ to $\cL_n$ with a large passage time or it starts from outside $S$, in which case $\Gamma$ must have large fluctuations. Let
        \begin{align*}
            \cD &= \bigcup_{x \in S}\bc{T^*_{x, n} \geq (4n + {\delta} t^2 2^{1/3} n^{1/3})\frac{K-k}{K}}&& \mathrm{and}&& \cD' = \bc{\Gamma \cap \cL_{k,K} \notin S}.
        \end{align*}
    The events $\cA_i$ are all increasing events, hence their intersection $\cA$ is also an increasing event. By definition, event $\cD$ is also an increasing event. Hence, by the BK inequality
        \begin{align*}
             \P \br{ \cA \:\square\: \cD} \leq \P \br{ \cA} \P \br{\cD}, 
        \end{align*}
    where $\cA \:\square\: \cD$ denotes the disjoint occurrence of $\cA$ and $\cD$. Furthermore, dividing the $2n^{2/3}$ length interval $S$ into intervals of size $(n(K-k)/K)^{2/3}$, using \cite[Proposition A.5]{BBF22} and taking a union bound we have
    \[
        \P \br{\cD} \leq (\frac{K}{K-k})^{2/3}e^{- \frac{ct^3 \delta^{3/2}(K-k)}{K}}.
    \]
    Consider the event
    \[
        \cE_0 =\bc{\abs{\psi(z^* - z^{k-1})} \geq \frac{n^{2/3}}{10}}.
    \]
    The point $z^*$ lies on the geodesic $\Gamma_{z^{k-1}, z^k}$. Thus on the event $\cE_0$, $\Gamma_{z^{k-1}, z^k}$ must have transversal fluctuations of the order $n^{2/3}$. However, these fluctuations are typically of the order $(n/K)^{2/3}$ and by a simple application of \cite[Proposition 2.1]{BBB23} (see also \cite[Theorem 3]{BSS17B}) we have that 
    \begin{equation} \label{eq: order}
        \P \br{\cE_0} \leq e^{-cK^2}.
    \end{equation}
    Now, consider the geodesics 
    \begin{align*}
        \Gamma^1 = \Gamma_{z^{k-1} + \frac{n^{2/3}}{5} 1_f,n}^* &&\mathrm{and}&& \Gamma^2= \Gamma_{z^{k-1} - \frac{n^{2/3}}{5} 1_f,n}^*.
    \end{align*}
    For $i =1,2$, the geodesic $\Gamma^i$ stays close to the line $\mathbb{L}_i = \bc{w \:\vert\: \psi(w) = \psi(z^{k-1}) -(-1)^i 2n^{2/3}/5}$. We define the maximal transversal fluctuation of $\Gamma^i$ before $\cL{k,K}$ as
    \[
        \mathrm{TF}^i := \sup \bc{\abs{\psi(w) - \psi(z^{k-1}) +(-1)^i 2n^{2/3}/5} \:\vert\: w \in \Gamma^i, \phi(w) \leq \frac{2nk}{K}}.
    \]
    Let $\cE_i$ be the event that $\bc{\mathrm{TF}^i > n^{2/3}/5}$. The quantity $\mathrm{TF}^i$ will typically be of the order $(n/K)^{2/3}$ as well. Hence by using arguments similar to the ones used for $\cE_0$, we can show that these events also have probability less than $e^{-cK^2}$. However, if none of the unlikely events $\cE_0, \cE_1, \cE_2$ happen then the paths $\Gamma^1, \Gamma$ and $\Gamma^2$ will be ordered from left to right and will cause the path $\Gamma$ to pass through $S$. Thus,
    \[
      \P \br{\cD'} \leq 3 e^{-cK^2}.
    \]
    Now, on the event $\cB_k^1 \cap \cA$ either the event $\cD \:\square\: \cA$ or $\cD'$ happens. Therefore,
        \begin{align} \label{eq: bulk} 
            \P \br{\cB_k^1 \:\vert\: \cA} \leq \frac{\P \br{\cD \:\square\: \cA} + \P \br{\cD'}}{\P\br{\cA}} \leq (\frac{K}{K-k})^{2/3}e^{- \frac{ct^3 \delta^{3/2}(K-k)}{K}} + \frac{e^{-cK^2}}{\P\br{\cA}}.
        \end{align}
    To bound the second event, note that
        \begin{align*}
            \cB_{k}^2 &\subseteq \bc{l(\Gamma_{z^{k-1},\cL_{k, K}}) \geq \frac{4n}{K} + \delta t^2 2^{4/3} n^{1/3}\frac{K-k}{2K}}.
        \end{align*}
        Observe that this latter event depends only on $\cA^1_{k}\cap \cA^2_k$ (if $k\le m$) or $\cA_{k}$ (otherwise) and are independent of other events defininig $\cA$.
    Let $a := \inf_i \bc{\P \br{\cA_i^2 \cap \cA_i^1}, \P(\cA_i)}$. We have, by \cite[Lemma A.4]{BBF22}
        \begin{equation} \label{eq: pathend}
        \begin{split}
            \P \br{\cB_{k}^2 \:\vert\: \cA} &\leq \P \br{l(\Gamma_{z^{k-1},\cL_{k, K}}) \geq \frac{4n}{K} + \delta t^2 2^{4/3} n^{1/3}\frac{K-k}{2K} \: \vert\: \cA}\\
            &\leq \frac{\P \br{T^*_{\frac{n}{K}} \geq \frac{4n}{K} + \delta t^2 2^{4/3} n^{1/3}\frac{K-k}{2K}}}{a} \\
             &\leq \frac{e^{-\frac{ct^3 \delta^{3/2}(K-k)^{3/2}}{K}}}{a} .
        \end{split}
        \end{equation}
    We shall apply a series of union bounds to get our result. Firstly, we have that $\P \br{\cB_k}$ is bounded above by the sum of three quantities,
    \begin{align*}
        \P\br{\cB'_{k} \:\vert\: \cA} \leq \frac{e^{-cK^2}}{\P(\cA)} + K^{2/3}e^{-\frac{ct^3 \delta^{3/2}(K-k)}{K}} + \frac{e^{- \frac{c t^3 \delta^{3/2}(K-k)^{3/2}}{K}}}{a} .
    \end{align*}    
    This implies that,
    \begin{align*}
        \P\br{\cup_{k \leq M}\cB'_{k} \:\vert\: \cA} \leq \sum_{k < M} \br{\frac{e^{-cK^2}}{\P(\cA)} + K^{2/3}e^{-\frac{ct^3 \delta^{3/2}(K-k)}{K}} + \frac{e^{- \frac{c t^3 \delta^{3/2}(K-k)^{3/2}}{K}}}{a}}.
    \end{align*}
    We shall individually bound all the three sums. The first term is simply upper bounded by
    \[
        \frac{Me^{-cK^2}}{\P(\cA)},
    \]
    which can be made arbitrarily small by taking $t$ large due to our choice of parameters in \eqref{para}. The next sum can also be made arbitrarily small since by our choice of $\delta$ and using $K-k\ge t$ each of the individual terms in the sum is upper bounded by
    $$K^{2/3}\exp(-ct^{0.9})$$
    and there are at most $O(t^{5/2})$ terms in the sum. 
  
   For the last sum, observe that by the FKG inequality,
   $$\P\br{\cA_i^1 \cap \cA_i^1} \geq \P \br{\cA_i^1} \P\br{\cA_i^1} \geq e^{- \frac{ct^3}{K}}.$$
   Thus, $a \geq e^{- \frac{ct^3}{K}}$. Using the choice of parameters we get that the numerator in each term is upper bounded by $\exp(-ct^{3.9})$. As before, the sum can be made arbitrarily small by choosing $t$ large.

    These inequalities show that  $\P\br{\cup_{k \leq M}\cB'_{k} \:\vert\: \cA}\le \frac{1}{8}$ if $t$ is sufficiently large. To bound the conditional probability of the event $\cB_{k}\setminus\cB'_{k}$ given $\cA$ for $k<m$, notice that in this case there exists a path $\Gamma$ as before which starts at $x^{k-1}_1$ instead. Arguing as before, we show  $\P\br{\cup_{k \leq m}(\cB_{k}\setminus \cB'_{k}) \:\vert\: \cA}\le \frac{1}{8}$ if $t$ is sufficiently large. Adding, we get $\P(\cB_{v_2}\mid \cA)\le 1/4$ as required. 
\end{proof}

To complete the proof we need to show that the probability that $\Gamma^*_{x^M,n}$ ends outside $I + v_{n,t}^*$ is small given $\cA$. This is the content of the following lemma. 
\begin{lemma} \label{l:epfluc}
    For the choice of parameters in \eqref{para}, we have the following bound for transversal fluctuations
    \[
        \P \br{\Gamma^*_{x^M, n}(n) \notin I + v_{n,t}^* \: \vert \: \cA} \leq \frac{1}{4}.
    \]
\end{lemma}
\begin{proof}
  Notice that the event in the statement of the lemma only depends on the passage times above the anti-diagonal line through $x^M$. Therefore, its conditional probability given $\cA$ is the same as its conditional probability given $\cA_{M+1}$. Thus
    \begin{equation*} 
        \begin{split}
            \P \br{\bc{\Gamma^*_{x^M, n}(n) \notin I + v_{n,t}^*} \mid \cA} &= \P \br{\bc{\Gamma^*_{x^M, n}(n) \notin I + v_{n,t}^*} \mid \cA_{M+1}} \\
            &\leq \dfrac{\P \br{T_{x^M, J} \geq \br{4n + \delta t^2 2^{4/3} n^{1/3}}\frac{K-M}{K}}}{\P(\cA_{M+1})},
        \end{split}
    \end{equation*}
where $J = \cL_n \setminus \bc{I + v_{n,t}^*}$. By \cite[Lemma A.3]{BBF22},
\begin{align*}
    \P \br{T_{x^M, J} \geq \br{4n + \delta t^2 2^{4/3} n^{1/3}}\frac{K-M}{K}} \leq e^{-c \frac{\delta^3 K^2}{(K-M)^2}}\le \exp(-ct^{1.8})
\end{align*} 
where the final step is obtained by using our choice of parameters. Observe also the we have from \eqref{e:alb} that $\P(\cA_{M+1})\ge \exp(-c\frac{(K-M)\delta^{3/2}t^3}{K})\ge \exp(-ct^{0.9})$. 

Thus we have that,
\begin{align*}
    \P \br{\Gamma^*_{x^M, n}(n) \notin I + v_{n,t}^* \: \vert \: \cA} &\leq \frac{\P \br{T_{x^M, J} \geq \br{4n + \delta t^2 2^{4/3} n^{1/3}}\frac{K-M}{K}}}{\P \br{\cA_{M+1}}}\\
    & \leq \frac{1}{4}
\end{align*}
where the final inequality is obtained by choosing $t$ sufficiently large. 
\end{proof}

\begin{proof}[Proof of Proposition \ref{p:p2llbbasic}]
     If both the geodesics $\Gamma^*_{n, v_1}$ and $\Gamma^*_{n, v_2}$ coalesce with $\Gamma^*_{n, x^M}$ and $\Gamma^*_{n, x^M}$ ends in the interval $I + v_{n,t}^*$ then the event $\cA_*$ will happen. By Lemma \ref{l:epfluc} and Lemma \ref{l:coal}, 
    \begin{align*}
        \P \br{\cA_* \:\vert \: \cA} &\ge 1 - \P\br{\Gamma^*_{n, v_1} \;\mathrm{does\;not\;intersect}\; \gamma \;\mathrm{after}\; \cL_{M,K} \:\vert \: \cA}\\
        &\;\;\;\;\;\;\; - \P\br{\Gamma^*_{n, v_2}\; \mathrm{does \;not \;intersect} \;\gamma\; \mathrm{after}\; \cL_{M,K} \:\vert \: \cA} \\
        &\;\;\;\;\;\;\; - \P \br{\Gamma^*_{x^M, n}(n) \notin I + v_{n,t}^* \: \vert \: \cA}\\
        &\ge 1 - \frac{3}{4}.
    \end{align*}
    This gives us that, 
    \[
        \P\br{\cA_*} \geq \frac{1}{4} \P \br{\cA} \ge \frac{1}{4} \exp\br{-\br{\frac{4}{3}+\epsilon}  t^3},
    \]
    completing the proof.
\end{proof}

\begin{proof}[Proof of Corollary \ref{c:p2l}]
    Corollary \ref{c:p2l} follows from Theorem \ref{t:p2l} by using an identical argument to the one used in Section \ref{s:p2p} to derive Corollary \ref{c:liu} from Theorem \ref{t:p2p}. We omit the details to avoid repititions. 
\end{proof}

\section{Point-to-point geodesics}
\label{s:p2p}
We shall use the same averaging argument as before. As explained in the introduction, we shall prove something stronger for the upper bound. Recall that for $\gamma\in (0,1)$, we defined the point $v_{n,\gamma,t}=(\gamma n-tn^{2/3}, \gamma n+tn^{2/3})$. Recall $I=I_{\delta}$ from the previous section. Let us denote $\widetilde{I}_{\gamma}$ to be the subset of points $v\in I$ such that the geodesic $\Gamma_{v,\mathbf{n}+v}$ passes through $v+v_{n,\gamma, t}$. Recall that $c_{\gamma}=3^{-1}(\gamma(1-\gamma))^{-3/2}$, and $c_{1/2}=8/3$. We have the following upper bound.

\begin{proposition}
    \label{p:p2pub}
    In the set-up of Theorem \ref{t:p2p}, for a fixed $\gamma\in (0,1)$ there exists $c_0>0$ such that for $c\in (0,c_0)$ and $\delta= \frac{c}{t}>0$ we have
    $$\E[|\widetilde{I}_{\gamma}|]\le  \exp(-(c_{\gamma}-\epsilon)t^3).$$
\end{proposition}

For the lower bound we have the following result in the case $\gamma=\frac{1}{2}$. 

\begin{proposition}
    \label{p:p2plb}
    In the set-up of Theorem \ref{t:p2p}, there exists $\delta>0$ such that
    $$\P[|\widetilde{I}_{1/2}|\ge 1]\ge  \exp(-(\frac{8}{3}+\epsilon)t^3).$$
\end{proposition}

Clearly Theorem \ref{t:p2p} follows from the case $\gamma=1/2$ of Proposition \ref{p:p2pub} and Proposition \ref{p:p2plb}.

\subsection{Upper bound}

Analogous to Lemma \ref{l:p2lub} we have the following lemma. 

\begin{lemma}
    \label{l:p2pub}
    In the set-up of Theorem \ref{t:p2p}, there exists $\delta=\frac{c}{t}>0$ such that
    $$\P[|\widetilde{I}_{\gamma}|\ge 1]\le  \exp(-(c_{\gamma}-\epsilon)t^3).$$
\end{lemma}

\begin{proof}
    As in the proof of Lemma \ref{l:p2lub}, using Proposition \ref{p:inf}, we conclude that for $C$ large enough 
    $$\P[\inf_{v\in I} T_{v,\mathbf{n}+v}\le 4n-Ctn^{1/3}]\le \exp(-c_{\gamma}t^{3}).$$
    Therefore it suffices to show that 
    $$\P[\sup_{v\in I} T_{v,v+v_{n,\gamma,t}}+ T_{v+v_{n,\gamma,t},v+\mathbf{n}}\ge 4n-Ctn^{1/3}]\le \exp(-(c_{\gamma}-\epsilon)t^{3}).$$
    Call this event $\cB$. 
Observe first that 
$$\sup_{v\in I} T_{v,v+v_{n,\gamma,t}}+ T_{v+v_{n,\gamma,t},\mathbf{n}}\le \sup_{v\in I,u\in v_{n,\gamma,t}+I} T_{v,u}+\sup_{u\in v_{n,\gamma,t}+I, w\in \mathbf{n}+I} T_{u,w}$$
and the two terms on the right hand side above are independent.\footnote{To make these terms independent we can just ignore the last vertex weight in a path, all the tail estimates for passage times still hold true.}
For $s>0$, let $A^1_{s}$ denote the event that 
$$\sup_{v\in I}T_{v,v+v_{n,\gamma,t}}\ge 4\gamma n -\frac{t^2}{(2\gamma)^{4/3}}2^{4/3}(\gamma n)^{1/3}+s 2^{4/3}n^{1/3}.$$
Similarly, let $A^2_{s}$ denote the event that 
$$\sup_{v\in I}T_{v+v_{n,\gamma,t}+v+\mathbf{n}}\ge 4(1-\gamma) n -\frac{t^2}{(2(1-\gamma))^{4/3}}2^{4/3}((1-\gamma) n)^{1/3}+s 2^{4/3}n^{1/3}.$$

It follows that 
\begin{equation}
    \label{e:sum1}\P(\cB) \le \sum_{s_1,s_2} \P(A_{s_1} \cap B_{s_2})
\end{equation}
where the sum is over all $s_1,s_2\in \N$ such that $s_1+s_2\ge \frac{t^2}{2^{4/3}\gamma(1-\gamma)}-Ct$. We shall sum over a larger (for $t$ large) space  
$s_1+s_2\ge \frac{t^2(1-\epsilon)^{2/3}}{2^{4/3}\gamma(1-\gamma)}$.
From Proposition \ref{p:sup} we have $\P(A^1_{s})\le \exp(-\frac{4(1-\epsilon)s^{3/2}}{3\gamma^{1/2}})$ and $\P(A^2_{s}) \le \exp(-\frac{4(1-\epsilon)s^{3/2}}{3(1-\gamma)^{1/2}})$. Therefore, using independence
$$\P(A^{1}_{s_1}\cap A^2_{s_2})\le \exp \left(-\frac{4(1-\epsilon)}{3}(s_1^{3/2}/\gamma^{1/2}+s_2^{3/2}/(1-\gamma)^{1/2})\right).$$
Next we try to find the largest term in \eqref{e:sum1}. We need to find 
$$\min  \frac{4(1-\epsilon)}{3}\left(\frac{s_1^{3/2}}{\gamma^{1/2}}+ \frac{s^{3/2}_2}{(1-\gamma)^{1/2}}\right)$$
subject to the constraints $s_1,s_2\ge 0$ and $s_1+s_2 \ge (1-\epsilon)^{2/3}\frac{t^2}{2^{4/3}\gamma(1-\gamma)}$.  It is easy to see that the minimum must occur on the line $s_1+s_2=(1-\epsilon)^{2/3}\frac{t^2}{2^{4/3}\gamma(1-\gamma)}$. Using the substitution $s_1=(1-\epsilon)^{2/3}at^{2}$ and $s_2=(1-\epsilon)^{2/3}bt^{2}$, it is equivalent to minimize 
$$\frac{a^{3/2}}{\gamma^{1/2}}+\frac{b^{3/2}}{(1-\gamma)^{1/2}}$$
subject to the constraint $a+b=\frac{1}{2^{4/3}\gamma(1-\gamma)}$. It can be checked (e.g.\ using Lagrange multipliers) that the minimum is attained at $a=\frac{1}{2^{4/3}(1-\gamma)}$ and $b=\frac{1}{2^{4/3}\gamma}$, and the minimum is $\frac{(1-\epsilon)^2 t^3}{3(\gamma(1-\gamma))^{3/2}}$. Since $\epsilon$ is arbitrary it follows that each term on the right hand side of \eqref{e:sum1} is upper bounded by $\exp(-(c_{\gamma}-\epsilon)t^3)$ and since the number of terms is $O(t^2)$ the lemma follows. 
\end{proof}

To complete the proof of Proposition \ref{p:p2pub} we need the following lemma. 

\begin{lemma}
    \label{l:multirarenew}
    There exists $k_0,c>0$ such that for $k_0\le k \le n^{0.01}$
    $$\P[|\widetilde{I}_{\gamma}|\ge k]\le \exp(-ck^{1/128}).$$
\end{lemma}

\begin{proof}
    Notice that if $v_1, v_2\in \widetilde{I}_{\gamma}$ it follows that either $\Gamma_{v_1,v_1+v_{n,\gamma,t}}$ and $\Gamma_{v_2,v_2+v_{n,\gamma,t}}$ are disjoint or $\Gamma_{v_1+v_{n,\gamma,t},\mathbf{n}}$ and $\Gamma_{v_2+v_{n,\gamma,t},\mathbf{n}}$ are disjoint. Therefore it follows, if $|\widetilde{I}_{\gamma}|\ge k$
    that there exists a set $\hat{I}\subset \widetilde{I}$ with $|\hat{I}|\ge \frac{k}{2}$ such that we have either for all $v\in \hat{I}$, $\Gamma_{v,v+v_{n,\gamma,t}}$ are pairwise disjoint or for all $v\in \hat{I}$, $\Gamma_{v+v_{n,\gamma,t}, \mathbf{n}}$ are pairwise disjoint. As in Lemma \ref{l:multirare}, using \cite[Proposition 1.6]{BBB23} we get $\P[|\widetilde{I}|\ge k]\le \exp(-ck^{1/128})$, as required. 
\end{proof}

\begin{proof}[Proof of Proposition \ref{p:p2pub}]
The proof is completed following an argument identical to the one in the proof of Proposition \ref{p:p2lub} using Lemmas \ref{l:p2pub} and \ref{l:multirarenew} instead of Lemmas \ref{l:p2lub} and \ref{l:multirare}; we omit the details.    
\end{proof}

\subsection{Lower bound}
For the lower bound, the construction and the proof techniques are similar to the ones in Section \ref{s:p2l}. Many events/parameters in the construction we are going to present will be analogous to events/parameters in the previous section. As such, we have used the same notations for them. Let $v_1$ and $v_2$ denote the left and right endpoint of $I$ respectively. Recall that $v_{n,t}=(\frac{n}{2}-tn^{2/3},\frac{n}{2}+tn^{2/3})$. Let $\cA_*$ denote the event that $\Gamma_{v_1, v_1+\mathbf{n}}$ intersects the line $\cL_{n/2}$ at a point below $v_1+v_{n,t}$ and $\Gamma_{v_2,\mathbf{n}+v_2}$ intersects $\cL_{n/2}$ at a point above $v_2+v_{n,t}$. By ordering of the geodesics, it is clear that, on $\cA_*$ we have $|\widetilde{I}_{1/2}|\ge 1$. Therefore Proposition \ref{p:p2plb} follows from the next result. 

\begin{proposition}
    \label{p:p2plbbasic}
    In the set-up of Proposition \ref{p:p2plb}, there is $\delta<1 $ such that $\P(\cA_*)\ge \exp(-(\frac{8}{3}+\epsilon)t^{3})$. 
\end{proposition}

\begin{figure}[!tbp]
    \centering
    \includegraphics[scale = 0.5]{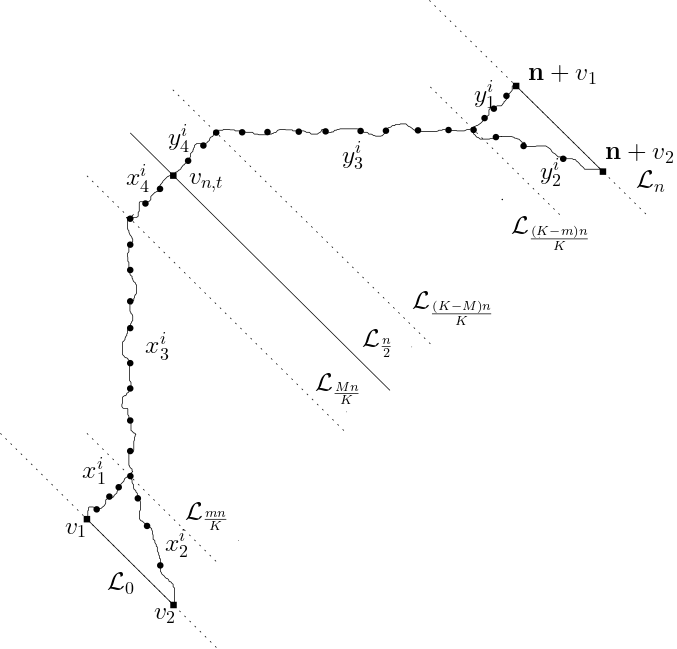}
    \caption{Illustration of the construction in the proof of Proposition \ref{p:p2plbbasic}: as before we construct two planted paths (that are the same except near endpoints) one each from $v_1$ and $v_2$ to $v_1+\mathbf{n}$ and $v_2+\mathbf{n}$ respectively such that the lengths of these paths are larger than $4n$, and that these paths roughly follow the piecewise linear paths joining $v_i$ to $v_i+v_{n,t}$ to $v_{i}+\mathbf{n}$. We further ask that the excess weight is approximately uniformly distributed along different parts of these paths. We then show that the probability cost of planting these paths is $\approx\exp(-8t^3/3)$ and that the event $\cA_*$ happens with probability bounded away from $0$ conditional on the planted paths.}
    \label{fig: p2p_exact}
\end{figure}

Our choice of parameters is similar to the one in Section \ref{s:p2l}.
\begin{align} \label{para-pp}
    \delta &= t^{-0.4}, &K=  t^{2.5},\\
    M &= \frac{K}{2} - t, & m = \frac{M}{t^{1.4}}.  
\end{align} 
For ease of notation, we define $\kappa = K/2M \approx 1$ for $t$ large. Each parameter in our construction has a similar role as before. We first build paths from $v_i$ to $\cL_{n/2}$ in a similar way as in Section \ref{s:p2l}. However, we shall exclusively use point to point geodesics to reach $\cL_{n/2}$ as opposed to the previous section which made use of both point to point and point to line geodesics. We ask that the geodesics connecting consecutive points have more than the typical weight. Let
\begin{align*}
    \cA^1_i &= \bc{T_{x^{i-1}_1,x^{i}_1} \geq 4n/K + 4 \delta t^2 n^{1/3}/K} &\mathrm{where}& &x^{i}_1 =  \frac{i\bfn}{K} + \delta n^{2/3} 1_f\\
    \cA^2_i &= \bc{T_{x^{i-1}_2,x^{i}_2} \geq 4n/K + 4 \delta t^2 n^{1/3}/K} &\mathrm{where}&&x^{i}_2 = \frac{i \bfn}{K} + \frac{(2i-m)}{m}\delta n^{2/3} 1_f,\\
    \cA^3_i &= \bc{T_{x_3^{i-1},x_3^{i}} \geq 4n/K + 4 \delta t^2 n^{1/3}/K} &\mathrm{where}&&x_3^{i} = \frac{i \bfn}{K} + \frac{i}{M}tn^{2/3}1_f,\\
    \cA^4_i &= \bc{T_{x_4^{i-1},x_4^{i}} \geq 4n/K + 4 \delta t^2 n^{1/3}/K} &\mathrm{where}&&x_4^{i} = \frac{i\bfn}{K} + tn^{2/3}1_f.
\end{align*}
Here, we define the events $\cA^1_i$ and $\cA^2_i$ for $0 < i \leq m$, events $\cA^3_i$ for $m < i \leq M$ and events $\cA^4_i$ for $M < i \leq K/2$. Next, we shall reflect all points $x_j^i$ used in the above definitions about the line $\cL_{n/2}$ and ask that equivalent events hold for their reflected counterparts. For $k= 1,2,3$ and $4$, let us denote by $y^i_k$, the reflection of $x^i_k$ by the line $\cL_{n/2}$. The ``reflected" event corresponding to $\cA_i^1$ would then be 
\[
     \bc{T_{y^{i}_1,y^{i-1}_1} \geq 4n/K + 4 \delta t^2 n^{1/3}/K}, 
\]
and we can define the reflections for the remaining events similarly. For $K/2<i \le K$, let us denote by ${\cA}^{k}_{i}$, the ``reflected" event corresponding to $\cA^{k}_{K+1-i}$, whenever the latter is well defined. The complete construction is illustrated in Figure \ref{fig: p2p_exact}.

We shall call the concatenation of all these paths as $\gamma$ and the intersection of all the constructed events $\cA$. The part of $\gamma$ from $v_i$ to $\bfn + v_i$ is denoted by $\gamma_i$ for $i=1,2.$ We now show that the probability cost for planting $\gamma$ as above is of the same order as the desired bound. By Proposition \ref{t:tails}, for arbitrarily small $\epsilon$, we can choose $t$ large enough such that
\begin{align*}
    \P\br{\cA^1_i} &\ge \exp\br{-\br{\frac{8}{3} + \frac{\epsilon}{2}} \frac{\delta^{3/2}t^3}{K}},
    & \P\br{\cA^2_i} \ge \exp\br{-\br{\frac{8}{3} + \frac{\epsilon}{2}} \frac{t^3 (\delta + 4\kappa ^2)^{3/2}}{K}},\\
    \P\br{\cA^3_i} &\ge \exp\br{-\br{\frac{8}{3} + \frac{\epsilon}{2}} \frac{t^3 (\delta+\kappa ^2)^{3/2}}{K}},
    &\P\br{\cA^4_i} \ge \exp \br{-\br{\frac{8}{3} + \frac{\epsilon}{2}} \frac{\delta^{3/2} t^3}{K}},
\end{align*}
for $i\le K/2$, and the reflected events ${A}^{k}_{i}$ for $i>K/2$ have identical probability lower bounds. Since all of these events are increasing we have by the FKG inequality that
\begin{align*}
    \P\br{\cA} &\geq \exp\br{-c\delta t^2}\exp\br{-\br{\frac{8}{3} + \frac{\epsilon}{2}} (\delta+\kappa ^2)^{3/2}t^3}\\
    & \geq \exp\br{-\br{\frac{8}{3} + {\epsilon}} t^3},    
\end{align*}
for $t$ large enough, depending on $\epsilon$. Since, we will be conditioning on these events throughout all our arguments, we adopt the following notation:
\[
    \bar{\P}\br{\cdot} = \P\br{\cdot \;\vert\;\cA}.
\]

We shall show that $\Gamma_{v_2, \bfn + v_2}$ lies to the left of $v_{n,t} + v_2$ with large probability when conditioned on $\cA$. To do so, we prove probability bounds on the last/first intersection point of $\Gamma_{v_2, \bfn + v_2}$ with $\gamma$ before/after $\cL_{n/2}$. As before, we perform the following relabeling:
\begin{equation*}
    z^i = \begin{cases}
        x^i_2 \;\;\;\;\;\;\;\;\;\;\;\; \mathrm{for} \;\;\;\;\;\;\;\; 0 \le i \leq m,\\
        x^i_3 \;\;\;\;\;\;\;\;\;\;\;\; \mathrm{for} \;\;\;\;\;\;\;\; m < i \leq M,\\
        x^i_4 \;\;\;\;\;\;\;\;\;\;\;\; \mathrm{for} \;\;\;\;\;\;\;\; M < i \leq K/2,\\
        y^{K-i}_4 \;\;\;\;\;\;\;\; \mathrm{for} \;\;\;\;\;\;\;\; K/2 < i \leq K-M,\\   
        y^{K-i}_3 \;\;\;\;\;\;\;\; \mathrm{for} \;\;\;\;\;\;\;\; K-M < i \leq K-m,\\
        y^{K-i}_2 \;\;\;\;\;\;\;\; \mathrm{for} \;\;\;\;\;\;\;\; K-m < i \leq K.\\
    \end{cases}
\end{equation*}
Let us say the last intersection point of $\Gamma_{v_2, \bfn + v_2}$ with $\gamma$ before $\cL_{n/2}$ was the point $x\in \gamma_2$ and let $k^*$ be the minimum number such that $x$ lies before $\cL_{k^*, K}$. Let the first intersection point of $\Gamma_{v_2, \bfn + v_2}$ with $\gamma$ after $\cL_{n/2}$ be $y\in \gamma_2$, correspondingly let $m^*$ be the largest number such that $\phi(y) \geq \phi(z^{m^*})$. Let $\cC_{k,k'}$ be the event that $k^* = k$ and $m^* = k'$.
 
On the event $\cC_{k,k'}$, we look at the path $\Gamma$ from $z^{k-1}$ to $z^{k'+1}$ which follows $\gamma_2$ from $z^{k-1}$ to $x$, $\Gamma_{v_2,\mathbf{n}+v_2}$ from $x$ to $y$ and $\gamma_2$ again from $y$ to $z^{k'+1}$. By definition, $\Gamma$ does not intersect  $\gamma$ between $z^k$ and $z^{k'}$ and has an excess length of at least $4 \delta t^2 n^{1/3} \frac{({k'} - k)}{K}$.

We divide $\Gamma$ into three parts $\Gamma_\mathrm{left}, \Gamma_\mathrm{mid}, \Gamma_\mathrm{right}$ such that 
    \begin{align*}
    \Gamma_\mathrm{left} &= \bc{\phi(z^{k-1}) \leq \phi(w) \leq \phi(z^k): w \in \Gamma},
    \\ \Gamma_\mathrm{mid} &= \bc{\phi(z^k) \leq \phi(w) \leq \phi(z^{k'}) : w \in \Gamma},
    \\ \Gamma_\mathrm{right} &= \bc{\phi(z^{k'}) \leq \phi(w) \leq \phi(z^{k'+1}): w \in \Gamma}.
    \end{align*}
    One of the following three events must happen:
    \begin{align*}
        \cC_{\rightarrow k} &= \bc{\ell(\Gamma_\mathrm{left}) \geq 2(\phi(z^k) - \phi(z^{k-1})) + 2\delta t^2 n^{1/3}\biggl(\frac{K/2-k}{K}\biggr)}, \\
        \cC_{k \rightarrow {k'}} &= \bc{\ell(\Gamma_\mathrm{mid}) \geq 2(\phi(z^{k'}) - \phi(z^k)) + 2\delta t^2 n^{1/3}\br{\frac{{k'} - k - \alpha}{K}}}, \\
        \cC_{{k'} \rightarrow } &= \bc{\ell(\Gamma_\mathrm{right}) \geq 2(\phi(z^{{k'}+1}) - \phi(z^{k'})) + 2\delta t^2 n^{1/3}\br{\frac{{k'} - K /2+ \alpha}{K}}},
    \end{align*}
    where $\alpha$ will be appropriately chosen to deal with the case when ${k' - K/2}$ is very small in the later arguments. For the first event, we have that
    \begin{align*}
        \P\br{\cC_{\rightarrow k}} &\le \P\br{{T}^*_{\bfn\frac{(k-1)}{K}, \frac{nk}{K}} \geq \frac{4n}{K} + 2\delta t^2 n^{1/3}\frac{K/2-k}{K}}.
    \end{align*}
    By equation \eqref{eq: pathend}, for $a = \inf_i \bc{\P\br{\cA^3_i}, \P\br{\cA^1_i \cap \cA^2_i}}$, we have that
    \begin{align} \label{eq: cC-k}
        \bar{\P}\br{\cC_{\rightarrow k}} \leq \frac{1}{a}C e^{- \frac{ct^3 \delta^{3/2}(K/2-k)^{3/2}}{K}}.
    \end{align}
    By symmetry we will have a similar bound for $\cC_{{k'} \rightarrow }$,
    \begin{align} \label{eq: cCm-}
        \bar{\P}\br{\cC_{{k'} \rightarrow }} \leq \frac{1}{a}C e^{-\frac{ct^3 \delta^{3/2}({k'}- K/2 + \alpha )^{3/2}}{K}}.
    \end{align}
    We define the following sets
    \begin{align*}
        S &= \bc{\abs{\psi(w-z^{k-1})} \leq n^{2/3} \:\vert\: w \in \cL_{k,K}} && \mathrm{and} && S' &= \bc{\abs{\psi(w-z^{{k'}+1})} \leq n^{2/3} \:\vert\: w \in \cL_{{k'},K}}.
    \end{align*}
    The event $\cC_{k \rightarrow {k'}}$ implies that either there is a path from $S$ to $S'$ disjoint from $\gamma$ with appropriate excess lengths or $\Gamma$ has large fluctuation near its endpoints. The chain of arguments is now the same as in the proof of Lemma \ref{l:coal}. Using the BK inequality, the $\bar{\P}$ probability of the first case can be upper bounded by the probability that there exists a large path across an $\frac{n(k'-k)}{K}\times n^{2/3}$ parallelogram. We bound this by dividing $S$ and $S'$ into intervals of length $(n(k'-k)/K)^{2/3}$, use \cite[Theorem 4.2]{BGZ21} and taking a union bound. For the second term note that since the restriction of $\Gamma$ between $x$ and $y$ is a geodesic. It follows that the typical fluctuation of $\Gamma$ at the lines $\cL_{k,K}$ and $\cL_{k',K}$ are of the order $(\frac{n}{K})^{2/3}$ around the point where the straight line joining $x$ and $y$ intersects these lines respectively. We now recall the proof of \eqref{eq: order}. Using the same ordering argument along with \cite[Proposition 2.1]{BBB23} (see also \cite[Theorem 3]{BSS17B}), we get that the $\bar{\P}$ probability $\Gamma$ that does not intersect $S$ or $S'$ is at most $\P(\cA)^{-1}e^{-cK^2}$.
    
    Combining these,
    \begin{align} \label{eq: cCk-{k'}}
        \bar{\P}\br{\cC_{k \rightarrow {k'}}} \leq C (\frac{K}{k'-k})^{4/3} e^{-  \frac{ct^3 \delta^{3/2}({k'} - k -\alpha)}{K}} + \frac{1}{\P(\cA)} e^{-cK^2}.
    \end{align}

\begin{lemma} \label{l: lip}
    For the choice of parameters as in \eqref{para-pp}, we have that for $t$ large enough
    \[
        \Pc{\bigcup_{k < M, k' \ge K/2} \cC_{k,k'}} \leq \frac{1}{10}.
    \]
\end{lemma}
\begin{proof}
    \begin{align*}
        \Pc{\bigcup_{k < M, {k'} \ge K/2} \cC_{k,{k'}}} &\leq \Pc{\bigcup_{k < M, {k'} \ge K/2} \cC_{\rightarrow k}\cup \cC_{k \rightarrow {k'}} \cup \cC_{{k'} \rightarrow }}\\
        &\leq \Pc{\bigcup_{k \leq M} \cC_{\rightarrow k}} + \Pc{\bigcup_{k < M, {k'} \ge K/2} \cC_{k \rightarrow {k'}}} + \Pc{\bigcup_{{k'} \ge K/2} \cC_{{k'} \rightarrow }} \\
        &\leq \sum_{k \leq M}\Pc{ \cC_{\rightarrow k}} + \sum_{k < M, {k'} \ge K/2}\Pc{\cC_{k \rightarrow {k'}}} + \sum_{{k'} \ge K/2}\Pc{\cC_{{k'} \rightarrow }}.
    \end{align*}
    The second inequality follows from realizing that $\cC_{\rightarrow k}$ depends only on $k$ and $\cC_{{k'} \rightarrow }$ depends only on ${k'}$. We shall now bound each of these three sums by 1/30. By equation \eqref{eq: cCk-{k'}},
    \begin{align*}
        \sum_{k < M, {k'} \ge K/2}\Pc{\cC_{k \rightarrow {k'}}} &\leq \sum_{k < M, {k'} \ge K/2} K^{4/3} e^{-\frac{ct^3 \delta^{3/2}({k'} -k- \alpha)}{K}} + \sum_{k < M, {k'} \ge K/2} \frac{1}{\P(\cA)} e^{-cK^2}.
    \end{align*}
    The above quantity can be made arbitrarily small by taking $t$ large enough and imposing the condition that  $\alpha$ is smaller than $\frac{K}{2}-M$ (by a factor) and the fact that $K$ is a polynomial in $t$. By  \eqref{eq: cC-k},
    \begin{align*}
        \sum_{k \leq M}\Pc{\cC_{\rightarrow k}} &\leq \sum_{k \leq M} \frac{1}{a}C e^{- \frac{ct^3 \delta^{3/2}(K/2-k)^{3/2}}{K}} \\
        &= e^{- \frac{ct^3 (\delta + 4\kappa ^2)^{3/2}}{K}} . \frac{e^{-  \frac{ct^3 \delta^{3/2}(K/2-M)^{3/2}}{K}}}{1- e^{-  \frac{ct^3 \delta^{3/2}(K/2-M)^{1/2}}{K}}}.
    \end{align*}
    To bound the right hand side by $1/30$, we will need $K/2-M$ to be much larger than $\delta^{-1}$ which can be done by taking $t$ to be large. Similarly, for $\cC_{{k'} \rightarrow }$ we require that $\alpha \gg \delta^{-1}$. We accomplish this by letting $\alpha = \frac{1}{2} \br{K/2 -M}$ and choosing $t$ appropriately large.
\end{proof}
The above lemma (essentially) tells us that the last intersection point of $\Gamma_{v_2,v_2+\mathbf{n}}$ with $\gamma$ before $\cL_{n/2}$ has a very small probability of being being before $z^M$. By symmetry, we have the same statement for the first intersection point of $\Gamma_{v_2,v_2+\mathbf{n}}$ with $\gamma$ after $\cL_{n/2}$. Combining these we have the following lemma.
\begin{lemma} \label{l: lfip}
    For the choice of parameters as in \eqref{para-pp}, we have that for $M' = K-M$ and $t$ large enough,
    \[
    \Pc{\bigcup_{k \geq M, k' \leq M'} \cC_{k,k'}} \geq \frac{3}{4}.
    \]
\end{lemma}
\begin{proof}
    From Lemma \ref{l: lip}, we have that
    \begin{equation} \label{eq: lip1}
        \Pc{\bigcup_{k < M, k' \geq K/2} \cC_{k,k'}} \leq \frac{1}{10}.
    \end{equation}
    We know that $\Gamma_{v_2,v_2+\mathbf{n}}$ must intersect $\gamma$ before and after $\cL_{n/2}$. However the complement of the event in \eqref{eq: lip1} does not imply that the last intersection of $\Gamma_{v_2,v_2+\mathbf{n}}$ and $\gamma$ before $\cL_{n/2}$ is after $\cL_{nM/K}$. This is because the last intersection of $\gamma$ with $\Gamma_{v_2,v_2+\mathbf{n}}$ before $\cL_{n/2}$ might belong to $\gamma_1$ and not $\gamma_2$ (notice that this can only happen in $k^*<m$ but this case was excluded in the definition of $\cC_{k,k'}$). However, the $\bar{\P}$-probability of this can be shown to be smaller than $\frac{1}{40}$ by the same argument in Lemma \ref{l: lip} with a modified choice of $z^i$; we omit the details. 
    It therefore follows that 
     \[
        \Pc{\bigcup_{k \ge M, k' \geq K/2} \cC_{k,k'}} \ge \frac{7}{8}.
    \]
    By symmetry, we also have the analogous statement for the first intersection point after $\cL_{n/2}$:
    \[
        \Pc{\bigcup_{k' \le M', k \leq K/2} \cC_{k,k'}} \ge \frac{7}{8}.
    \]
The claimed result follows from a union bound. 
\end{proof}
What we need to show now is that on the event described in Lemma \ref{l: lfip}, it is likely that $\Gamma_{v_2,v_{2+\mathbf{n}}}$ passes to the left of $v_2+v_{n,t}$. To do so, we shall shift our focus to a nearby geodesic. We adopt the symbol $\gamma_\mathrm{top}$ to denote the segments of $\gamma$ which lie between $x^M_4$ and its reflection (equal to $x^{M'}_4$) about $\cL_{n/2}$; see Figure \ref{fig:p2p_top}. 
\[
    \gamma_\mathrm{top} = \bc{w : w \in \gamma,\: \phi(x^M_4) \leq \phi(w) \leq \phi(x^{M'}_4)}.
\]
\begin{proof}[Proof of Proposition \ref{p:p2plbbasic}]
    We draw lines parallel to $x+y=0$ from the points $x^{M}_4$ and $x^{M'}_4$. We are interested in the intersection points of these lines with the line $y-x = (2t- \delta)n^{2/3}$, call them $x'$ and $y'$. Therefore,
    \[
    \phi(y')-\phi(x') = \phi(x^{M'}_4)-\phi(x^{M}_4) = \frac{4n(M'-M)}{K}.
    \]
    By transversal fluctuation estimates (see e.g.\ \cite[Proposition C.9]{BGZ21}),
    \begin{align*}
        \P\br{2t n^{2/3}\ge \abs{\Gamma_{x',y'}} \ge (2t -2\delta)n^{2/3}} 
        & \geq 1 - C e^\frac{-c K^2 \delta^3}{(M'-M)^2}.
    \end{align*}
   We denote the event in the L.H.S. by $\cD$. By our choice of parameters in \eqref{para-pp},
   \begin{align*}
       \Pc{\cD^c} \leq \frac{\P\br{\cD^c}}{\P\br{\bigcap_{M<i<M'}\cA_i^4}} \leq  \frac{ Ce^\frac{-c K^2 \delta^3}{(M'-M)^2}}{e^{-\frac{8}{3} \frac{(M'-M)\delta^{3/2} t^3}{K}}} \leq 1/8.
   \end{align*} 
    Call the event that $\Gamma_{v_2,\bfn+v_2}$ intersects $\gamma_\mathrm{top}$ at least once before $\cL_{n/2}$ and after it as $\cC_{v_2}$. By Lemma \ref{l: lfip}, we have that
    \[
        \Pc{\cC_{v_2}} \geq  \Pc{\bigcup_{k \geq M, k' \leq M'} \cC_{k,k'}}  \geq 3/4.
    \]
    \begin{figure}
        \centering
        \includegraphics[scale = 0.4]{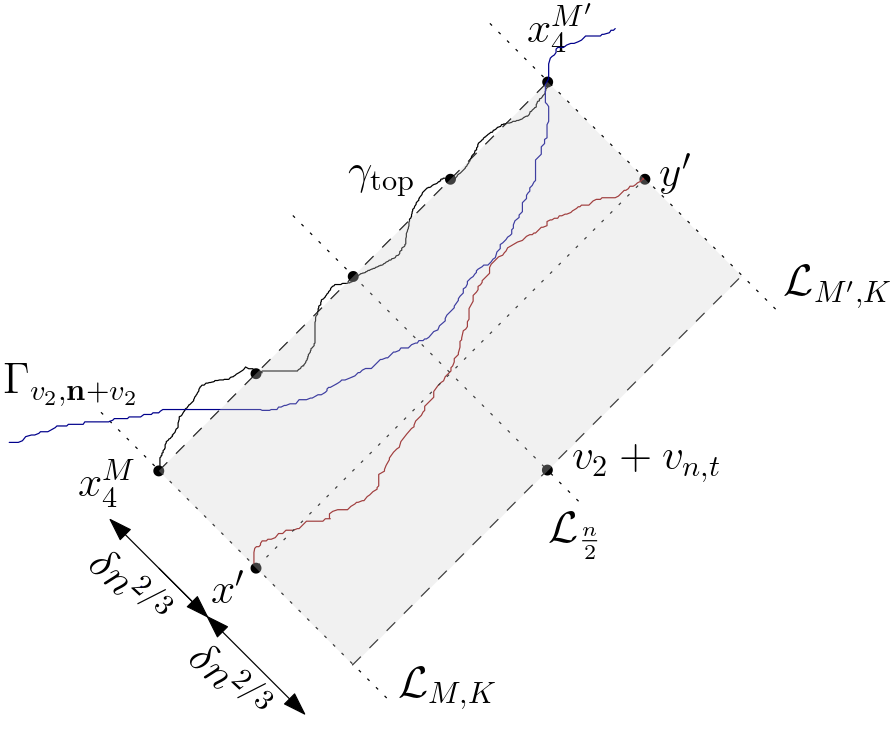}
        \caption{The red path represents $\Gamma_{x',y'}$ and the blue path represents $\Gamma_{v_2,\bfn+v_2}$. On the event $\cD$, $\Gamma_{x',y'}$ will be forced to be in the shaded region which implies, by ordering of geodesics, that if $\Gamma_{v_2,v_2+\mathbf{n}}$ intersects $\gamma_{\mathrm{top}}$ once before and once after $\cL_{n/2}$ then it must be above the red path at $\cL_{n/2}$.}
        \label{fig:p2p_top}
    \end{figure}
    By ordering of the geodesics, on the event $\cD$ all the geodesic segments present in $\gamma_\mathrm{top}$ shall lie to the left of $\Gamma_{x',y'}$ (refer to Figure \ref{fig:p2p_top}). Thus on the event $\cC_{v_2} \cap \cD$, $\Gamma_{v_2,\bfn+v_2}$ is forced to be to the left of $\Gamma_{x',y'}$ on an interval around $\cL_{n/2}$. Since $\Gamma_{x',y'}$ intersects $\cL_{n/2}$ to the left of the point $v_2+v_{n,\gamma, t}$, we have that
    \[
        \Pc{\abs{\Gamma_{v_2,\bfn+v_2}(n/2)} \geq \phi(v_2+v_{n,\gamma, t})} \geq \Pc{\cC_{v_2} \cap \cD} \geq 5/8.
    \]
    By our construction, using similar arguments for $\Gamma_{v_1,v_1+\mathbf{n}}$ we also get,
    \[
        \Pc{\abs{\Gamma_{v_1,\bfn+v_1}(n/2)} \leq \phi(v_1+v_{n,\gamma, t})} \geq 5/8.
    \]
    Thus,
    \[
    \Pc{\cA_*} \geq \Pc{\bc{\abs{\Gamma_{v_1,\bfn+v_1}(n/2)}\leq \phi(v_1+v_{n,\gamma, t})} \cap \bc{\abs{\Gamma_{v_2,\bfn+v_2}(n/2)}\geq \phi(v_2+v_{n,\gamma, t})}} \geq \frac{1}{4}.
    \]
    This implies that,
    \[
    \P \br{\cA_*} \geq \Pc{\cA^*} \P\br {\cA} \geq  \frac{1}{4}\exp\br{-\br{\frac{8}{3} + {\epsilon}} t^3};
    \]
    completing the proof. 
\end{proof}

We now present a short proof for Corollary \ref{c:liu}.
\begin{proof}[Proof of Corollary \ref{c:liu}]
    The lower bound follows directly from Proposition \ref{p:p2plb}. For the upper bound we shall employ a union bound. Let $N \in \mathbb{N}$ and consider the interval, $$J_N = \bc{(-x, x) : 0 < x < N t n^{2/3}} \subset \cL_0.$$
    By taking a union bound and applying Theorem \ref{t:p2p}, we have that for any $\epsilon >0$ and $n,t$ sufficiently large
    \[
    \P\br{\Gamma_n(n/2) \in J_N + v_{n,t}} \leq N t \exp\br{-\br{\frac{8}{3}-\frac{\epsilon}{2}}t^3}.
    \]
    Next, by applying \cite[Proposition 2.1]{BBB23} (see also \cite[Theorem 3]{BSS17B}), there exists an absolute constant $c>0$ such that
    \[
    \P\br{\Gamma_n(n/2) \text{ lies above } v_{n,{(N+1)t}}} \leq \exp \br{-cN^3 t^3}.
    \]
    Now by taking $N$ large enough, we can claim that
    \[
    \P\br{\Gamma_n(n/2) \text{ lies above } v_{n,{(N+1)t}}} \leq \exp \br{-\br{\frac{8}{3} -\frac{\epsilon}{2}} t^3}.
    \]
    Thus, we have that
    \begin{align*}  
        \P\br{\Gamma_n(n/2) \text{ lies above } v_{n,{t}}} &\leq \P\br{\Gamma_n(n/2) \text{ lies above } v_{n,{(N+1)t}}} + \P\br{\Gamma_n(n/2) \in J_N + v_{n,{t}}} \\
        &\leq \exp \br{-\br{\frac{8}{3} -\epsilon} t^3}.
    \end{align*}
    Where the last inequality follows from taking $n,t$ large enough, depending on $\epsilon$.
\end{proof}

 \bibliography{bib}

\begin{thebibliography}{10}

\bibitem{Aga23}
Pranay Agarwal.
\newblock Lower bound for large local transversal fluctuations of geodesics in last passage percolation, 2023.
\newblock arXiv 2311.00360.

\bibitem{BLS12}
Jinho Baik, Karl Liechty, and Gr{\'e}gory Schehr.
\newblock On the joint distribution of the maximum and its position of the {A}iry$_2$ process minus a parabola.
\newblock {\em Journal of mathematical physics}, 53(8), 2012.

\bibitem{BBB23}
M.~Bal{\'a}zs, R.~Basu, and S.~Bhattacharjee.
\newblock Geodesic trees in last passage percolation and some related problems.
\newblock {\em \url{https://arxiv.org/abs/2308.07312}}, 2023.

\bibitem{BBBK24}
Jnaneshwar Baslingker, Riddhipratim Basu, Sudeshna Bhattacharjee, and Manjunath Krishnapur.
\newblock Optimal tail estimates in $\beta$-ensembles and applications to last passage percolation, 2024.
\newblock arXiv 2405.12215.

\bibitem{BBF22}
R.~Basu, O.~Busani, and P.~Ferrari.
\newblock On the exponent governing the correlation decay of the $\text{Airy}_1$ process.
\newblock {\em Comm. Math. Phys.}, 398:1171--1211, 2022.

\bibitem{BGZ21}
R.~Basu, S.~Ganguly, and L.~Zhang.
\newblock Temporal correlation in last passage percolation with flat initial condition via brownian comparison.
\newblock {\em Comm. Math. Phys.}, 383:1805--1888, 2021.

\bibitem{BB21}
Riddhipratim Basu and Manan Bhatia.
\newblock Small deviation estimates and small ball probabilities for geodesics in last passage percolation.
\newblock {\em arXiv preprint arXiv:2101.01717}, 2021.

\bibitem{BHS18}
Riddhipratim Basu, Christopher Hoffman, and Allan Sly.
\newblock Nonexistence of bigeodesics in planar exponential last passage percolation.
\newblock {\em Communications in Mathematical Physics}, 389, 2022.

\bibitem{BSS17B}
Riddhipratim Basu, Sourav Sarkar, and Allan Sly.
\newblock Coalescence of geodesics in exactly solvable models of last passage percolation.
\newblock {\em Journal of Mathematical Physics}, 60(9):093301, 2019.

\bibitem{BSS14}
Riddhipratim Basu, Vladas Sidoravicius, and Allan Sly.
\newblock Last passage percolation with a defect line and the solution of the {S}low {B}ond {P}roblem.
\newblock Preprint arXiv 1408.3464, 2014.

\bibitem{BF08}
Alexei Borodin and Patrik~L. Ferrari.
\newblock Large time asymptotics of growth models on space-like paths {I}: Push{ASEP}.
\newblock {\em Electronic Journal of Probability}, 13:1380--1418, 2007.

\bibitem{DDV24}
Sayan Das, Duncan Dauvergne, and Vir\'ag Balint.
\newblock Upper tail large deviations of the directed landscape, 2024.
\newblock arXiv 2405.14924.

\bibitem{DOV18}
Duncan Dauvergne, Janosch Ortmann, and B{\'a}lint Vir{\'a}g.
\newblock The directed landscape.
\newblock {\em Acta. Math.}, 229, 2022.

\bibitem{DV21}
Duncan Dauvergne and Bálint Virág.
\newblock The scaling limit of the longest increasing subsequence, 2021.
\newblock arXiv 2104.08210.

\bibitem{HS18}
Alan Hammond and Sourav Sarkar.
\newblock Modulus of continuity for polymer fluctuations and weight profiles in poissonian last passage percolation.
\newblock {\em Electron. J. Probab.}, 25:38 pp., 2020.

\bibitem{LR10}
M.~Ledoux and B.~Rider.
\newblock Small deviations for beta ensembles.
\newblock {\em Electron. J. Probab.}, 15:1319--1343, 2010.

\bibitem{Liu22}
Zhipeng Liu.
\newblock One-point distribution of the geodesic in directed last passage percolation.
\newblock {\em Probability Theory and Related Fields}, 184(1):425--491, 2022.

\bibitem{Sch12}
Gr{\'e}gory Schehr.
\newblock Extremes of n vicious walkers for large n: Application to the directed polymer and kpz interfaces.
\newblock {\em Journal of Statistical Physics}, 149(3):385--410, 2012.

\bibitem{Shen24}
Xiao Shen.
\newblock Lower bound for large midpoint transversal fluctuations in the corner growth model, 2024.
\newblock arXiv:2402.16332.

\end{thebibliography}
 \bibliographystyle{plain}

\end{document}